\documentclass[preprint]{elsarticle}
\usepackage{graphicx}
\usepackage{amsfonts,amssymb,amsbsy,amsmath,amsthm}
\usepackage{ifthen}
\newtheorem{thm}{Teorem}

\begin{document}
\title{Model Order Reduction for Nonlinear Schr\"odinger Equation}


\author[rtv]{B.~Karas\"{o}zen,\corref{cor1}}
\ead{bulent@metu.edu.tr}

\author[focal]{C.~Akkoyunlu}
\ead{c.kaya@iku.edu.tr}

\author[els]{M.~Uzunca}
\ead{uzunca@metu.edu.tr}

\cortext[cor1]{Corresponding author}

\address[rtv]{Department of Mathematics and Institute of Applied
  Mathematics, Middle East Technical University,  06800 Ankara,
  Turkey}
\address[focal]{Department of Mathematics and Computer Sciences,
       Istanbul K\"{u}lt\"{u}r University, 34156, Istanbul, Turkey}
\address[els]	{ Department of Mathematics, Middle East Technical University,  06800 Ankara,
  Turkey}

\begin{abstract}
We apply the proper orthogonal decomposition (POD) to the nonlinear Schr\"{o}\-dinger (NLS) equation to derive a reduced order model. The NLS equation is discretized in space by finite differences and is solved in time by structure preserving symplectic  mid-point rule. A priori error estimates are derived for the POD reduced dynamical system. Numerical results for one and two dimensional NLS equations, coupled NLS equation with soliton solutions show that the low-dimensional approximations  obtained  by POD reproduce very well  the characteristic dynamics of the system, such as preservation of energy and the solutions.
\end{abstract}

\begin{keyword}
Nonlinear Schr\"{o}dinger equation; proper orthogonal decomposition; model order reduction; error analysis
\end{keyword}
\maketitle

\section{Introduction}

The nonlinear Schr\"{o}dinger (NLS) equation   arises as the model
equation with second order dispersion and cubic nonlinearity
describing the dynamics of slowly varying wave packets in
nonlinear fiber optics, in water waves and  in Bose-Einstein condensate theory.
We consider the NLS equation
 \begin{eqnarray} \label{nlsdenk}
 \mathrm{i}\Psi_t + \Psi_{xx} + \gamma \mid \Psi \mid^2\Psi = 0
  \end{eqnarray}
with the periodic boundary conditions $\Psi(x+L,t) = \Psi(x,t)$.  Here $\Psi=\Psi(x,t)$ is a complex valued function, $\gamma$ is a parameter and $\mathrm{i}=\sqrt{-1}$.  The NLS equation is called "focusing" if $\gamma >0$ and "defocusing" if $\gamma<0$; for $\gamma=0$, it reduces to the linear Schr\"{o}dinger equation.   In last two decades, various numerical methods were applied for solving NLS equation, among them are the well-known symplectic and multisymplectic integrators and discontinuous Galerkin methods.

There is a strong need for model order reduction techniques to reduce the
computational costs and storage requirements in large scale simulations, yielding low-dimensional
approximations for the full high-dimensional dynamical system, which reproduce
the characteristic dynamics of the system. Among the model order reduction techniques, the proper orthogonal decomposition (POD) is one of the most widely used method.
It was first introduced for analyzing cohorent structures and turbulent flow in numerical simulation of fluid dynamics equations \cite{bertcs96}.
It has been successfully used in different fields including signal processing, fluid dynamics, parameter estimation, control theory and optimal control of partial differential equations.  In this paper, we apply the POD to the NLS equation. To the best of our knowledge, there is only one paper where POD is applied to NLS equation  \cite{schlizerman12}, where only one and two modes approximations of the NLS equation are used in the Fourier domain in connection with mode-locking ultra short laser applications. In this paper, the NLS equation being a semi-linear partial differential equation (PDE) is discretized in space and time by preserving the symplectic structure and the energy (Hamiltonian). Then, from the snapshots of the fully discretized dynamical system, the POD basis functions are computed using the singular value decomposition (SVD).  The reduced model consists of Hamiltonian ordinary differential equations (ODEs), which indicates that the geometric structure of the original system is preserved for the reduced model. The semi-disretized NLS equations and the reduced equations are solved in time using Strang splitting and mid-point rule. A priori error estimates are derived for POD reduced model, which is solved by mid-point rule.  It turns out that most of the energy of the system can be accurately approximated by
using few POD modes. Numerical results for a NLS equation with soliton solutions confirm that the energy of the system is well preserved by  POD approximation and the solution of the reduced model are close to the solution of the fully discretized system.\\
The paper is organized as follows. In Section 2, the  POD method and its application to semi-linear dynamical systems are reviewed.	In Section 3, a priori error estimators are derived for the mid-point time-discretization of semi-linear PDEs.  Numerical solution of the semi-discrete NLS equation and the POD reduced form are described in Section \ref{numnls}. In the last section, Section \ref{numres}, the numerical results for the reduced order models of NLS equations are presented.
\section{The POD approximation for semi-linear PDEs} \label{podover}
In the following, we briefly describe the important features of the POD reduced order modeling (ROM);
more details can be found in  \cite{kunisch01gpo}.
In the first step of the POD based model order reduction, the set of snapshots,
 the discrete solutions of the nonlinear PDE, are collected. The snapshots are usually equally spaced in time
 corresponding to the solution of PDE obtained by finite difference or finite element method.
  The snapshots are then used to determine the POD bases which are much smaller than the snapshot set.
  In the last step, the POD reduced order model is constructed to obtain approximate solutions of the PDE.
  We mention that the choice of the snapshots representing the dynamics of the underlying PDE is crucial for the effectiveness of POD based reduced model.\\
		Let $X$ be a real Hilbert space endowed with inner product
 $\left\langle \cdot,\cdot \right\rangle_X $ and norm $\left\|\cdot\right\|_X$. For $y_1,\ldots ,y_n\in X $, we set
$
V=span\left\{y_1,\cdot \cdot \cdot,y_n\right\},
$
 as the ensemble consisting of the snapshots $\left\{y_j\right\}_{j=1}^n$.
In the finite difference context, the snapshots can be viewed as discrete solutions  $y_j\in {\mathbb R}^m$ at time instances $t_j$, $j=1,\ldots ,n$, and $[y_1,\ldots ,y_n]\in {\mathbb R}^{m\times n}$ denotes the snapshot matrix.\\
Let $\left\{\psi_k\right\}_{k=1}^d$ denote an orthonormal basis of $V$ of dimension $d$. Then, any $y_j\in V$ can be expressed as
\begin{eqnarray}\label{poddis}
y_j=\sum_{k=1}^d\left\langle y_j,\psi_k\right\rangle_X\psi_k ,\quad j=1,\ldots ,n.
\end{eqnarray}
The POD is constructed by choosing the orthonormal basis such that for every $l\in\left\{1,\ldots ,d\right\}$, the mean square error between the elements $y_j$, $1\leq j\leq n$, and the corresponding $l-th$ partial sum of (\ref{poddis}) is minimized on average:
\begin{equation}\label{podbaz}
\min_{\tilde{u}_1,\ldots ,\tilde{u}_l \in X}\sum_{j=1}^n \alpha_j\left\|y_j-\sum_{k=1}^l\left\langle y_j,\tilde{u}_k\right\rangle_X\tilde{u}_k\right\|_X^2, \quad
\left\langle \tilde{u}_i,\tilde{u}_j\right\rangle_X=\delta_{ij},\quad 1\leq i,j\leq l .
\end{equation}
where $\alpha_j$'s are non-negative weights. Throughout this paper, we take the space $X=\mathbb{R}^m$ endowed with the weighted inner product $\langle u,v\rangle_W = u^TWv $ with the diagonal elements of the diagonal matrix $W$, and also $\alpha_j$'s are the trapezoidal weights so that we obtain all the computations in $L_2$-sense. Under these choices, the solution of the above minimization problem is given by the following theorem:
\begin{thm}
\cite{kunisch01gpo}. Let $Y=[y_1,\ldots ,y_n]\in\mathbb{R}^{m\times n}$ be a given matrix with rank $d\leq min\left\{m,n\right\}$. Further, let $Y=U\Sigma V^T$ be the SVD of $Y$, where $U=[u_1,\ldots ,u_m]\in\mathbb{R}^{m\times m}, V=[v_1,\ldots ,v_n]\in \mathbb{R}^{n\times n} $ are orthogonal matrices and the matrix $\Sigma\in \mathbb{R}^{m\times n}$ is all zero but first $d$ diagonal elements are the nonzero singular values, $\sigma_1\geq\sigma_2\geq\ldots \geq\sigma_d$, of Y. Then, for any $l\in\{ 1,\ldots ,d\}$, the solution to
\begin{equation}\label{thrm}
\min_{\tilde{u}_1,\ldots ,\tilde{u}_l \in\mathbb{R}^m}\sum_{j=1}^n \alpha_j\left\|y_j-\sum_{k=1}^l\left\langle y_j,\tilde{u}_k\right\rangle_{W}\tilde{u}_k\right\|_{W}^2, \quad
\left\langle \tilde{u}_i,\tilde{u}_j\right\rangle_{W}=\delta_{ij},\quad 1\leq i,j\leq l .
\end{equation}
is given by the singular vectors $\left\{u_i\right\}_{i=1}^l$.
\end{thm}
We consider the following initial value problem for POD approximation
\begin{equation}
\label{tsyst}
\dot{y}(t)=Ay(t)+f(t,y(t)),\quad t\in[0,T],\quad
y(0)=y_0,
\end{equation}
where $f:[0,T]\times\mathbb{R}^{m}\rightarrow \mathbb{R}^{m}$ is continuous in both arguments and locally Lipschitz-continuous with respect to the second argument. The semi-discrete form of NLS equation (\ref{nlsdenk}) is a semi-linear equation as (\ref{tsyst}) where the cubic nonlinear part is locally Lipschitz continuous.
Suppose that we have determined a POD basis $\left\{ \psi_j\right\}_{j=1}^l$ of rank $l\in \left\{ 1,\ldots ,d\right\}$ in $\mathbb{R}^{m}$, then we make the ansatz
\begin{eqnarray}\label{rom}
y^l(t)=\sum_{j=1}^l\underbrace{\left\langle y^l(t),\psi_j\right\rangle_W}_{=:\mathrm{y}_j^l(t)}\psi_j,\quad t\in[0,T].
\end{eqnarray}
Substituting (\ref{rom}) in (\ref{tsyst}), we obtain the reduced model
\begin{equation}\label{pod}
\sum_{j=1}^l\dot{\mathrm{y}}_j^l(t)\psi_j=\sum_{j=1}^l\mathrm{y}_j^l(t)A\psi_j+f(t,y^l(t)), \quad t\in[0,T], \quad
\sum_{j=1}^l\mathrm{y}_j^l(0)\psi_j=y_0.
\end{equation}
The POD approximation (\ref{pod}) holds after projection on the $l$ dimensional subspace $V^l=span\{\psi_1,\ldots ,\psi_l\}$. From (\ref{pod}) and $\left\langle \psi_j,\psi_i\right\rangle_W=\delta_{ij}$, we get 
\begin{eqnarray}\label{podd}
\dot{\mathrm{y}}_i^l(t)=\sum_{j=1}^l\mathrm{y}_j^l(t)\left\langle A\psi_j,\psi_i\right\rangle_W+\left\langle f(t,y^l(t)),\psi_i\right\rangle_W
\end{eqnarray}
for $1\leq i\leq l$ and $t\in(0,T]$. Let us introduce the matrix
$$
B=\{ b_{ij}\}\in \mathbb{R}^{l\times l}, \quad b_{ij}=\left\langle A\psi_j,\psi_i\right\rangle_W
$$
the non-linearity $F=(F_1,\cdots,F_l)^T:[0,T]\times \mathbb{R}^{l} \rightarrow \mathbb{R}^{l}$ by
$$
F_i(t,\mathrm{y})=\left\langle f(t,\sum_{j=1}^l\mathrm{y}_j\psi_j),\psi_i\right\rangle_W, \quad t\in[0,T], \quad \mathrm{y}=(\mathrm{y}_1,\cdots,\mathrm{y}_l)\in \mathbb{R}^{l}
$$
and the vector $\mathrm{y}^l=(\mathrm{y}_1^l,\ldots , \mathrm{y}_l^l)^T:[0,T] \rightarrow \mathbb{R}^{l}$. Then, (\ref{podd}) can be expressed as
\begin{equation}
\label{reduced}
\dot{\mathrm{y}}^l(t)=B\mathrm{y}^l(t)+F(t,\mathrm{y}^l(t)), \quad t\in(0,T].
\end{equation}
The initial condition of the reduced system is given by
$
\mathrm{y}^l(0)=\mathrm{y}_0
$
with
$$
\mathrm{y}_0=\left( \left\langle y_0,\psi_1\right\rangle_W, \ldots , \left\langle y_0,\psi_l\right\rangle_W \right)^T\in \mathbb{R}^{l}.
$$
The system  (\ref{reduced})  is called the POD-Galerkin projection for (\ref{tsyst}). The ROM is constructed with POD basis vectors $\{\psi_i\}_{i=1}^l$ of rank $l$. In case of $l<< d$, the $l-$dimensional reduced system
(\ref{reduced}) is a low-dimensional approximation for (\ref{tsyst}).\\
The POD basis can also be computed using eigenvalues and eigenvectors. We prefer singular value decomposition, because it is more accurate than the computation
of the eigenvalues. The singular values decay up to machine precision, where the eigenvalues stagnate several orders above due the fact
$\lambda_i =\sigma_i^2$ \cite{studinger13napod}. We notice that all singular values of the snapshot matrix $Y$ are normalized, so that
$\sum _{i=1}^m \sigma_i^2 =1$ holds.  .
The choice of $l$ is based on heuristic considerations combined with observing the ratio of the modeled energy to the total energy contained in the system $Y$ which is expressed by
the relative information content (RIC)
$$
\mathcal{E}(l)= \sum_{i=1}^l\sigma_i^2 \cdot 100 \%  .
$$
The total energy of the system is contained in a small number of  POD modes. In practice, $l$ is chosen by guaranteeing that $\mathcal{E}(l)$
capturing at least \% 99 of total energy of the system.
\section{POD error analysis for the mid-point rule}
A priori error estimates for POD method are obtained for linear and semi-linear parabolic equations in \cite{kunisch01gpo}, where the nonlinear part is assumed to be locally Lipschitz continuous as for the NLS equation. The error estimates derived for the backward Euler and Crank-Nicholson (trapezoidal rule) time discretization show that the error bounds depend on the number of POD basis functions.
Here, we derive the error estimates for the mid-point rule.  We apply the implicit midpoint rule for solving the reduced model (\ref{reduced}). By $Y_j$, we denote an approximation for $y^l$ at the time $t_j$. Then, the discrete system for the sequence $\{ Y_j\}_{j=1}^n$ in $V_n^l=\text{span}\{\psi_1^n,...,\psi_l^n\}$ ($l\leq d$) looks like
\begin{eqnarray}
\left\langle\frac{Y_j-Y_{j-1}}{\Delta t},\psi_i^n\right\rangle_W &=&\left\langle\frac{1}{2}A(Y_j+Y_{j-1})+f(t,\frac{Y_j+Y_{j-1}}{2}),\psi_i^n\right\rangle_W\\
\langle Y_1,\psi_i^n\rangle_W&=&\langle y_0,\psi_i^n\rangle_W,\qquad \qquad i=1,...,l
\end{eqnarray}
We are interested in estimating $\sum_{j=1}^n\alpha_j\left\|y(t_j)-Y_j\right\|_W^2$. For $u\in\mathbb{R}^m$, let us introduce the projection $P_n^l:\mathbb{R}^m\rightarrow V_n^l$ by
$$
P_n^lu=\sum_{i=1}^l\langle u,\psi_i^n\rangle_W\psi_i^n, \qquad \left\|P_n^l\right\|_{W}=1.
$$
We shall make use of the decomposition
$$
y(t_j)-Y_j=y(t_j)-P_n^ly(t_j)+P_n^ly(t_j)-Y_j=\varrho_j^l+\vartheta_j^l
$$
where $\varrho_j^l=y(t_j)-P_n^ly(t_j)$ and $\vartheta_j^l=P_n^ly(t_j)-Y_j$. Using that $\{\psi_i^n\}_{i=1}^l$ is the POD basis of rank $l$, we have the estimate for the terms involving $\varrho_j^l$
\begin{equation} \label{varte}
\sum_{j=1}^n\alpha_j\left\|y(t_j)-\sum_{i=1}^l<y(t_j),\psi_i^n>\psi_i^n\right\|_W^2=\sum_{j=1}^n\alpha_j\left\|y(t_j)-P_n^ly(t_j)\right\|_W^2=\sum_{j=1}^n\alpha_j
\left\|\varrho_j^l\right\|_W^2=\sum_{i=l+1}^d\sigma_i^2.
\end{equation}
Next,  we estimate the terms involving $\vartheta_j^l$. Using the notation $\bar{\partial}\vartheta_j^l=(\vartheta_j^l-\vartheta_{j-1}^l)/\Delta t$, we obtain
\begin{eqnarray}\label {dis}
\langle\bar{\partial}\vartheta_j^l,\psi_i^n\rangle_W &=&\left\langle P_n^l\left(\frac{y(t_j)-y(t_{j-1})}{\Delta t}\right)-\frac{Y_j-Y_{j-1}}{\Delta t},\psi_i^n\right\rangle_W \nonumber \\
&=&\left\langle \dot{y}(t_j-\frac{\Delta t}{2})-\left(\frac{1}{2}A(Y_j+Y_{j-1})+f(\frac{Y_j+Y_{j-1}}{2})\right),\psi_i^n\right\rangle_W \nonumber\\
&+&\left\langle P_n^l\left(\frac{y(t_j)-y(t_{j-1})}{\Delta t}\right)-\dot{y}(t_j-\frac{\Delta t}{2}),\psi_i^n\right\rangle_W  \\
&=&\left\langle A\left(y(t_j-\frac{\Delta t}{2})- \frac{Y_j+Y_{j-1}}{2}  \right)+f(y(t_j-\frac{\Delta t}{2}))-f(\frac{Y_j+Y_{j-1}}{2})+w_j^l+z_j^l,\psi_i^n\right\rangle_W
\nonumber
\end{eqnarray}
where
$$
z_j^l= P_n^l\left(\frac{y(t_j)-y(t_{j-1})}{\Delta t}\right)-\frac{y(t_j)-y(t_{j-1})}{\Delta t},\quad  w_j^l=\frac{y(t_j)-y(t_{j-1})}{\Delta t}-\dot{y}(t_j-\frac{\Delta t}{2})).
$$
Choosing $\psi_i^n=\vartheta_j^l+\vartheta_{j-1}^l$ in (\ref{dis}), we arrive at
\begin{eqnarray}\label{dis1}
\langle\bar{\partial}\vartheta_j^l,\vartheta_j^l+\vartheta_{j-1}^l\rangle_W &=&\langle A\left(y(t_j-\frac{\Delta t}{2})- \frac{Y_j+Y_{j-1}}{2}\right)+f(y(t_j-\frac{\Delta t}{2}))-f(\frac{Y_j+Y_{j-1}}{2}) \nonumber \\
& & +w_j^l+z_j^l,\vartheta_j^l+\vartheta_{j-1}^l\rangle_W.
\end{eqnarray}
Noting that
$$
\langle\bar{\partial}\vartheta_j^l,\vartheta_j^l+\vartheta_{j-1}^l\rangle_W =\frac{1}{\Delta t}\left(\left\|\vartheta_j^l\right\|_W^2-\left\|\vartheta_{j-1}^l\right\|_W^2\right)
$$
and using Lipschitz-continuity of $f$ and the Cauchy-Schwartz inequality in (\ref{dis1}), we get
\begin{equation} \label{vart}
\left\|\vartheta_j^l\right\|_W \leq \left\|\vartheta_{j-1}^l\right\|_W+\Delta t\left( (\left\|A\right\|_W+L_f) \left\|y(t_j-\frac{\Delta t}{2})-\frac{Y_j+Y_{j-1}}{2}\right\|_W+ \left\|z_j^l \right\|_W+ \left\|w_j^l \right\|_W\right).
\end{equation}
By Taylor series expansion
$$
y(t_j-\frac{\Delta t}{2}) = \frac{y(t_{j})+y(t_{j-1})}{2}+\frac{\Delta t}{2}\left( \dot{y}(\xi_{j-1})-\dot{y}(\xi_{j})\right)
$$
for some $\xi_j\in (t_j-\frac{\Delta t}{2},t_j)$ and $\xi_{j-1}\in (t_{j-1},t_j-\frac{\Delta t}{2})$. Then, we get
\begin{eqnarray} \label{taylor}
\left\|y(t_j-\frac{\Delta t}{2})-\frac{Y_j+Y_{j-1}}{2}\right\|_W &\leq& \frac{1}{2}\left( \left\|\varrho_j^l\right\|_W+\left\|\vartheta_j^l\right\|_W+\left\|\varrho_{j-1}^l\right\|_W+\left\|\vartheta_{j-1}^l\right\|_W+c_0\Delta t\right)
\end{eqnarray}
with $c_0=\dot{y}(\xi_{j-1})-\dot{y}(\xi_{j})$. Inserting (\ref{taylor}) in (\ref{vart}) and collecting the common terms yields
\begin{equation} \label{errmain}
(1-c_1\Delta t)\|\vartheta_j^l\|_W \leq (1+c_1\Delta t)\|\vartheta_{j-1}^l\|_W+\Delta t \left( c_1(\left\|\varrho_j^l\right\|_W+\left\|\varrho_{j-1}^l\right\|_W)+c_2\Delta t+\| z_j^l\|_W+ \| w_j^l\|_W\right)
\end{equation}
with $c_1=\text{max}\{ \| A\|_W, L_f\}$, $c_2= c_0c_1$. Moreover, for $0 < \Delta t \leq \frac{1}{2c_1}$, we have
$$
\frac{1}{1-c_1\Delta t} \leq 1+2c_1\Delta t
$$
and using the fact that $\Delta tj\leq T$, we get
\begin{equation} \label{deltat}
(1+2c_1\Delta t)^j \leq e^{2c_1T} , \qquad (1+c_1\Delta t)^j \leq e^{c_1T}.
\end{equation}
Summation on $j$ in (\ref{errmain}) by using (\ref{deltat}) and Cauchy-Schwarz inequality yields,
\begin{eqnarray} \label{varrho1}
\left\|\vartheta_j^l\right\|_W^2 &\leq& C\Delta t^2\sum_{k=1}^j\left( \left\|\varrho_k^l\right\|_W^2+\left\|\varrho_{k-1}^l\right\|_W^2+\Delta t^2+\| z_k^l\|_W^2+ \| w_k^l\|_W^2\right)
\end{eqnarray}
with $C=5e^{4c_1T}\text{max}\{ c_1^2,c_2^2,1,j\}$.
Next, we estimate the term involving $w_k^l$:
\begin{eqnarray} \label{estw}
\Delta t^2\sum_{k=1}^j\left\|w_k^l \right\|_W^2 &=& \Delta t^2\sum_{k=1}^j\left\| \frac{y(t_k)-y(t_{k-1})}{\Delta t}-\dot{y}(t_k-\frac{\Delta t}{2})\right\|_W^2 \nonumber \\[0.3cm]
&\leq& \tilde{C}^2\Delta t^4 \int_{0}^{T}\| y_{ttt}(t)\|^2dt
\end{eqnarray}
for a constant $\tilde{C}$ depending on $y$, but independent of $n$.
Now, we estimate the term involving $z_k^l$:
\begin{eqnarray} \label{estz}
\left\| z_k^l \right\|_W^2 &=& \left\| P_n^l\left( \frac{y(t_k)-y(t_{k-1})}{\Delta t}\right)-\frac{y(t_k)-y(t_{k-1})}{\Delta t} \right\|_W^2 \nonumber \\[0.3cm]
&\leq & 2\| w_k^l\|_W^2+4\| P_n^l\dot{y}(t_k)-\dot{y}(t_k)\|_W^2+4\left\| \dot{y}(t_k)-\frac{y(t_k)-y(t_{k-1})}{\Delta t}\right\|_W^2 \nonumber \\[0.3cm]
&\leq & 2\| w_k^l\|_W^2+4\| P_n^l\dot{y}(t_k)-\dot{y}(t_k)\|_W^2+4\left\| \dot{y}(t_k-\frac{\Delta t}{2})-\frac{y(t_k)-y(t_{k-1})}{\Delta t}\right\|_W^2+C_{Tay}\Delta t^2 \nonumber \\[0.3cm]
&\leq& 4\| P_n^l\dot{y}(t_k)-\dot{y}(t_k)\|_W^2+6\| w_k^l\|_W^2+C_{Tay}\Delta t^2
\end{eqnarray}
where $C_{Tay}=\ddot{y}(\xi)$ for some $\xi\in (t_k-\frac{\Delta t}{2},t_k)$.\\
For a sufficiently small $\Delta t$ satisfying $\Delta t\leq 2\alpha_k$ for $1\leq k\leq n$, we have
\begin{equation} \label{deltat2}
\Delta t^2\leq 2\alpha_k, \quad \Delta t^4\leq 2\alpha_k, \quad  \Delta t^2\sum_{k=1}^j\Delta t^2 \leq \sum_{k=1}^n2\alpha_k.
\end{equation}
Using (\ref{deltat2}) combining with (\ref{estw}) and (\ref{estz}), we arrive at
\begin{equation} \label{estwz}
\Delta t^2\sum_{k=1}^j\left( \| z_k^l \|_W^2 +\| w_k^l \|_W^2\right) \leq 8\sum_{k=1}^n\alpha_k\| P_n^l\dot{y}(t_k)-\dot{y}(t_k)\|_W^2+\hat{C}\Delta t^4
\end{equation}
with $\hat{C}=\text{max}\{ 14\tilde{C}^2 \| y_{ttt}\|_{L^2(0,T;\mathbb{R}^m)},2C_{Tay}\}$.
Imposing the estimates (\ref{estwz}) and (\ref{deltat2}) in (\ref{varrho1}), we obtain
\begin{equation} \label{varrho2}
\left\|\vartheta_j^l\right\|_W^2 \leq 4C\sum_{k=1}^n\alpha_k \|\varrho_k^l\|_W^2+C\sum_{k=1}^n\left( 2\alpha_k+8\alpha_k\| P_n^l\dot{y}(t_k)-\dot{y}(t_k)\|_W^2\right)+C\hat{C}\Delta t^4.
\end{equation}
In addition, we have that $\sum_{k=1}^n\alpha_k=T$ and $\| P_n^l\dot{y}(t_j)-\dot{y}(t_j)\|_W^2=\sum_{i=l+1}^d|\left\langle \dot{y}(t_j),\psi_i^n\right\rangle_W |^2$. Using these identities, we arrive at the estimate to the term involving $\vartheta_j^l$ as
\begin{equation} \label{varrho3}
\sum_{j=1}^n\alpha_j\|\vartheta_j^l\|_W^2 \leq C^*\left( \Delta t^4+\sum_{i=l+1}^d \left( \sigma_i^2+\sum_{j=1}^n\alpha_j|\left\langle \dot{y}(t_j),\psi_i^n\right\rangle_W |^2\right)\right)
\end{equation}
where $C^*=4CT\text{max}\{2T,8,\hat{C} \}$ and is dependent on $y$, $T$, but independent of $n$ and $l$.\\
Now, combining the estimates (\ref{varte}) and (\ref{varrho3}), we obtain finally the error estimate
\begin{eqnarray*}
\sum_{j=1}^n\alpha_j\| y(t_j)-Y_j\|_W^2 &=& \sum_{j=1}^n\alpha_j\| \vartheta_j^l+\varrho_j^l\|_W^2 \; \leq \; 2\sum_{j=1}^n\alpha_j\| \vartheta_j^l\|_W^2+ 2\sum_{j=1}^n\alpha_j\| \varrho_j^l\|_W^2\\[0.3cm]
&\leq& 2C^*\left( \Delta t^4+\sum_{i=l+1}^d \left( \sigma_i^2+\sum_{j=1}^n\alpha_j|\left\langle \dot{y}(t_j),\psi_i^n\right\rangle_W |^2\right)\right)+2\sum_{i=l+1}^d\sigma_i^2\\[0.3cm]
&\leq& C_E\left( \sum_{i=l+1}^d \left( 2\sigma_i^2+\sum_{j=1}^n\alpha_j|\left\langle \dot{y}(t_j),\psi_i^n\right\rangle_W |^2\right)+ \Delta t^4\right)
\end{eqnarray*}
where $C_E=\text{max}\{ 2,2C^*\}$ and is dependent on $y$, $T$, but independent of $n$ and $l$. As for the backward Euler and Cranck-Nicholson method
\cite{kunisch01gpo}, the error between the reduced and the unreduced solutions depend for the mid-point rule
on the time discretization and on the number of not modelled POD snapshots.
\section{Discretization of NLS equation} \label{numnls}
One dimensional NLS equation (\ref{nlsdenk}) can be written  by decomposing $\Psi=p+\mathrm{i}q$  in real and imaginary components
\begin{equation}\label{nlsorn1}
p_t   =   -q_{xx} -\gamma(p^2+q^2)q , \quad
q_t   =   p_{xx} + \gamma(p^2+q^2)p
\end{equation}
as an infinite dimensional Hamiltonian PDE in the phase space $u = (p,q)^T$
$$
\dot{u} = \mathcal{D} \frac{\delta\mathcal{H}}{\delta u}, \quad
\mathcal{H} = \int  \frac{1}{2}\left (p_x^2 +  q_x^2 - \frac{\gamma}{2} (p^2+q^2)^2
\right)dx,  \quad
\mathcal{D} = \left(  \begin{array}{cc} 0& 1  \\ -1 & 0 \end{array} \right).
$$
After discretizing the Hamiltonian in space by finite differences
\begin{equation} \label{hamilton}
H= \Delta x \sum_{j=1}^n \frac{1}{2} \left( \left( \frac{p_{j+1}-p_j}{\Delta x}\right)^2 +  \left( \frac{q_{j+1}-q_j}{\Delta x}\right)^2  - \frac{\gamma}{2} (p_j^2 +q_j^2)^2 \right).
\end{equation}
we obtain the semi-discretized  Hamiltonian ode's
\begin{equation} \label{nls}
p_t    =     -Aq - \gamma q(p^2+q^2), \quad
q_t    =   Ap+\gamma p(p^2+q^2),
\end{equation}
where $A$ is the circulant matrix
$$\label{B_matrisi}
A = \left(\begin{array}{ccccc}
-2 & 1 &  &  & 1 \\
1 & -2 & 1 &  &  \\
& \ddots & \ddots & \ddots &  \\
 &  & 1 & -2 & 1 \\
 1 &  &  & 1 & -2 \\
 \end{array}\right).
$$
\subsection{Reduced order model for NLS equation}
Suppose that we have determined POD bases $\left\{\psi_j\right\}_{j=1}^l$ and $\left\{\phi_j\right\}_{j=1}^l$ of rank $l=\left\{1,\ldots ,d\right\}$ in $\mathbb{R}^{m}$. Then we make the ansatz
\begin{eqnarray}\label{nls1d}
p^l=\sum_{j=1}^l\mathrm{p}_j(t) \psi_j(x),\quad q^l=\sum_{j=1}^l\mathrm{q}_j(t)\phi_j(x)
\end{eqnarray}
where $\mathrm{p}_j=\langle p^l,\psi_j\rangle_W,\quad \mathrm{q}_j=\langle q^l,\phi_j\rangle_W$. Inserting (\ref{nls1d}) into (\ref{nls}), and using
the orthogonality of the POD bases $\left\{\psi_j\right\}_{j=1}^l$ and $\left\{\phi_j\right\}_{j=1}^l$,  we obtain for $i=1,\cdots,l$ the systems
\begin{eqnarray*}
 \dot{\mathrm{p}_i}&=&-\sum_{j=1}^l\mathrm{q}_j\left\langle A\phi_j,\psi_i\right\rangle_W-\gamma\left\langle \left(\sum_{j=1}^l\mathrm{q}_j \phi_j\right)\left( \sum_{j=1}^l\mathrm{p}_j \psi_j\right)^2,\psi_i\right\rangle_W-\gamma\left\langle \left(\sum_{j=1}^l\mathrm{q}_j \phi_j\right)^3,\psi_i \right\rangle_W \\
 \dot{\mathrm{q}_i}&=&\sum_{j=1}^l\mathrm{p}_j\left\langle A \psi_j,\phi_i\right\rangle_W+\gamma\left\langle \left(\sum_{j=1}^l\mathrm{p}_j \psi_j\right)\left( \sum_{j=1}^l\mathrm{q}_j\phi_j\right)^2,\phi_i\right\rangle_W+\gamma\left\langle \left( \sum_{j=1}^l\mathrm{p}_j\psi_j\right)^3,\phi_i\right\rangle_W .
 \end{eqnarray*}
After defining $\Phi=[\phi_1,\phi_2,\cdots,\phi_l]\in\mathbb{R}^{m\times l},\; \Psi=[\psi_1,\psi_2,\cdots,\psi_l]\in\mathbb{R}^{m\times l},\; (B)_{ij} =\left\langle A\phi_j,\psi_i\right\rangle_W$, we obtain
\begin{equation}\label{dismor}
\begin{array}{lll}
 \dot{\mathrm{p}}&=&-B\mathrm{q}-\gamma \Psi^T\left( (\Phi\mathrm{q})\cdot (\Psi\mathrm{p})^2\right)-\gamma \Psi^T\left( (\Phi\mathrm{q})^3\right) \\
 \dot{\mathrm{q}}&=&B^T\mathrm{p}+\gamma \Phi^T\left( (\Psi\mathrm{p})\cdot (\Phi\mathrm{q})^2\right)+\gamma  \Phi^T\left( (\Psi\mathrm{p})^3\right)
 \end{array}
\end{equation}
with both '$\cdot$'  operation and the powers are hold elementwise.
The reduced NLS equation (\ref{dismor}) is also Hamiltonian and is solved, as the unreduced semi-discretized NLS equation (\ref{nlsdenk}), with  the symplectic midpoint method applying linear-nonlinear Strang splitting \cite{weidmann86sms}:
In order to solve (\ref{nls}) efficiently, we apply the  second order linear, non-linear Strang splitting
\cite{weidmann86sms}
$$\mathrm{i}u_t = {\mathcal N} u + {\mathcal L } u,\quad
{\mathcal L }u =  -u_{xx}, \quad {\mathcal N} u = -\gamma | u |^2 u.
$$
 The nonlinear parts  of the equations are solved by Newton-Raphson method. In the numerical examples, the boundary conditions are periodic, so that the resulting discretized matrices are circulant.   For solving the linear system of equations, we have used the Matlab toolbox {\bf smt} \cite{rediv012smt}, which is designed for solving linear systems with a structured coefficient matrix like the circulant and Toepltiz matrices. It reduces the number of floating point operations for matrix factorization to $O\;(n\log n)$.
\section{Numerical results}
\label{numres} All weights in the POD approximation are taken equally as $\alpha_i = 1/n$ and $W=I$. Then
the average ROM error, difference between the numerical solutions of NLS equation and ROM is measured in the form of the error between the fully discrete NLS solution
$$
\hbox{ROM error }  = \left(\frac{1}{n}\sum_{j=1}^n  \mid\mid y_h(t_j)-y_l(t_j)\mid\mid\right )^{1/2}.
$$
The average Hamiltonian ROM error is given by
$$
\left(\frac{1}{n}\sum_{j=1}^n  (H_h(t_j)-H_l(t_j)) ^2\right )^{1/2}.
$$
where $H_h(t_j)$ and $H_l(t_j)$ refer to the discrete Hamiltonian errors at the time instance $t_j$ corresponding to the full-order and ROM solutions, respectively. The energy of the Hamiltonian PDEs is usually expressed by the Hamiltonian. It is well known that symplectic integrators like the midpoint rule can preserve the only quadratic Hamiltonians exactly. Higher order polynomials and nonlinear Hamiltonians are preserved by the symplectic integration approximately, i.e. the approximate Hamiltonians do not show any drift  in long term integration.
For large matrices, the SVD is very time consuming. Recently several randomized methods are developed  \cite{halko11apc}, which are very efficient
when the rank is very small, i.e,  $d<< \min (m,n)$. We compare the efficiency of MATLAB programs {\it svd} and {\it fsvd} (based on the algorithm in  \cite{halko11apc}) for computation of singular values for
the NLS equations in this section, on a PC with AMD FX(tm)-8150 Eight-Core Processor and 32Gb RAM. The accuracy of the SVD is measured by $L_2$ norm,  $|| Y-U\Sigma V^T||_W$. The randomized version of SVD, the fast SVD {\it fsvd}, requires the rank of the matrix as input parameter, which can be determined by MATLAB's {\it rank} routine. When the singular values decay rapidly and the size of the matrices is very large, then randomized methods  \cite{halko11apc} are more efficient than MATLAB's {\bf svd}. Computation of the rank with {\it rank} and singular values with {\it fsvd} requires much less time than the {\it svd} for one and two dimensional NLS equations (Table 1).
\begin{table}[htb]
\begin{tabular}{|l|r|r|r|r|r|r|r|}
\hline
Problem & size of the matrix & rank & {\it rank} & {\it fsvd } & accuracy & {\it svd} & accuracy\\
\hline
1D NLS & 32 x 50001 & 15 & 0.14 & 0.73 & 6.4e-14 & 194.48 & 1.3e-16 \\
2D NLS & 6400 x 30001 & 25 & 279.74 & 3.81 & 2.01e-13 & 1300.42 & 2.47e-15 \\
CNLS & 128 x 2001 & 122 & 0.05 & 0.10 & 2.7e-14 & 1.19 & 1.5e-16 \\
\hline
\end{tabular}
\caption{Comparison of {\it svd} and {\it fsvd}}
\end{table}
\subsection{One-dimensional  NLS equation}
\begin{figure}[htb!]
\centering
\includegraphics[width=0.40\textwidth]{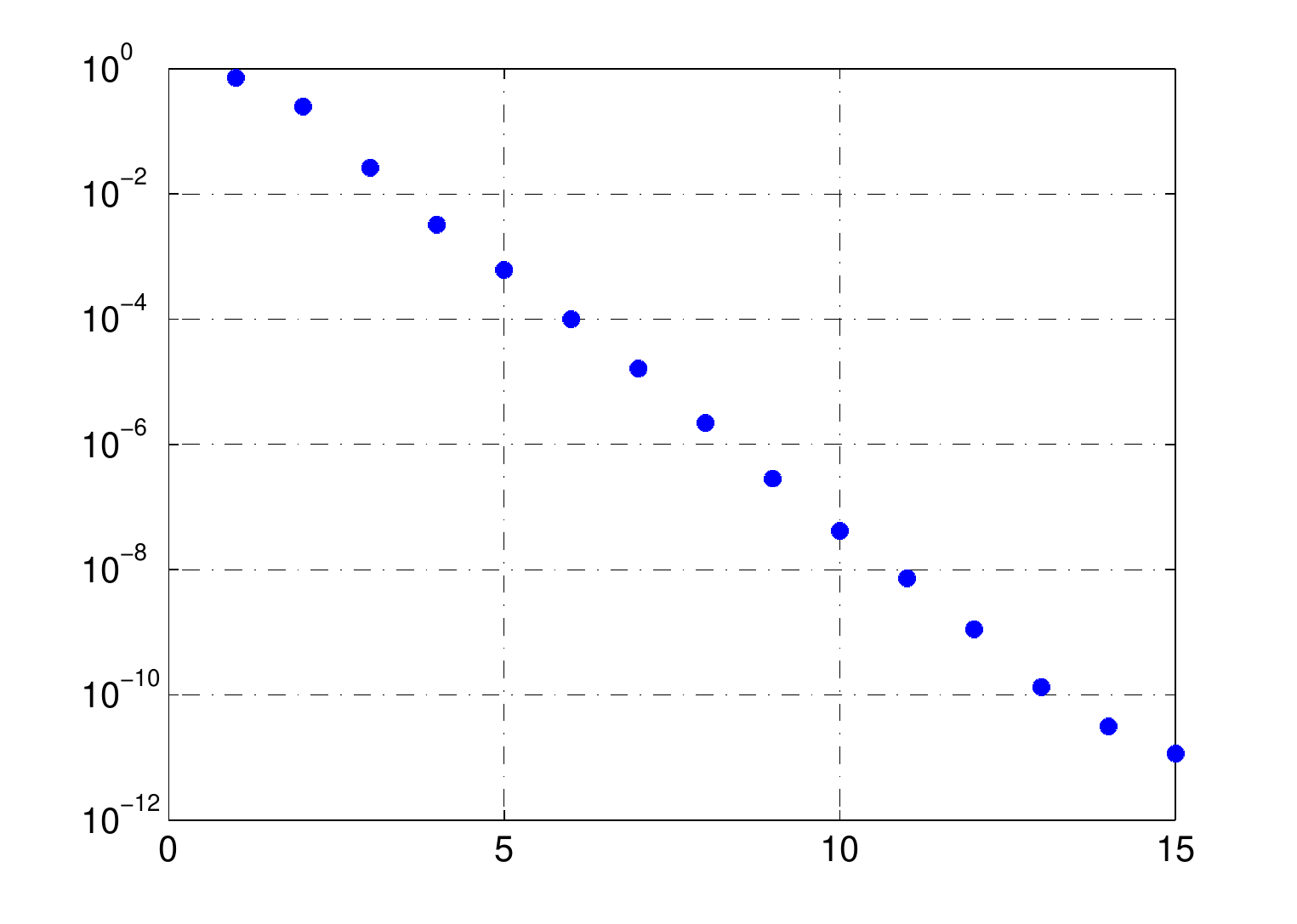}
\caption{1D NLS, Decay of the singular values }
\label{svd}
\end{figure}

\begin{figure}[htb!]
\centering
\includegraphics[width=0.40\textwidth]{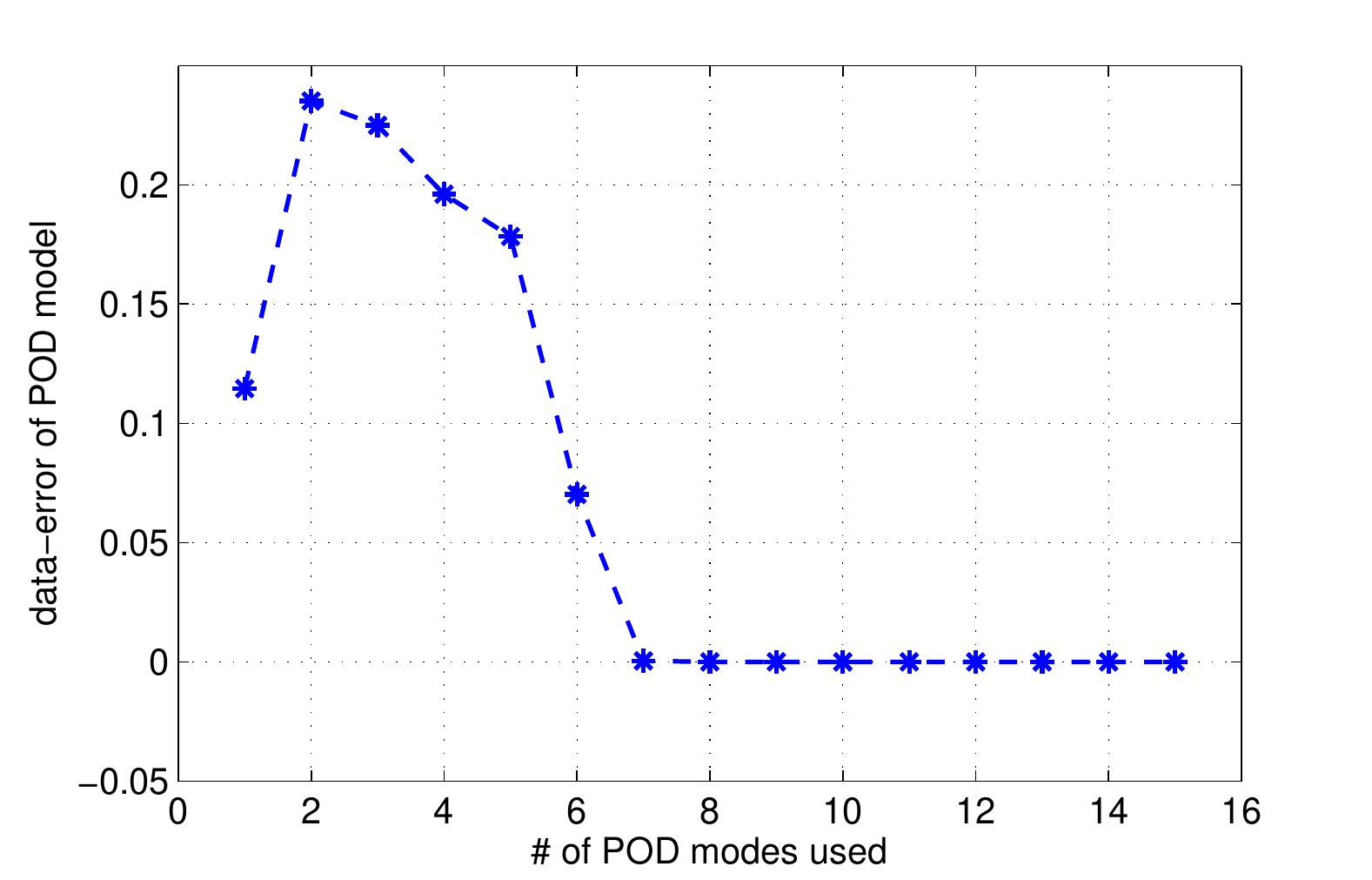}
\includegraphics[width=0.40\textwidth]{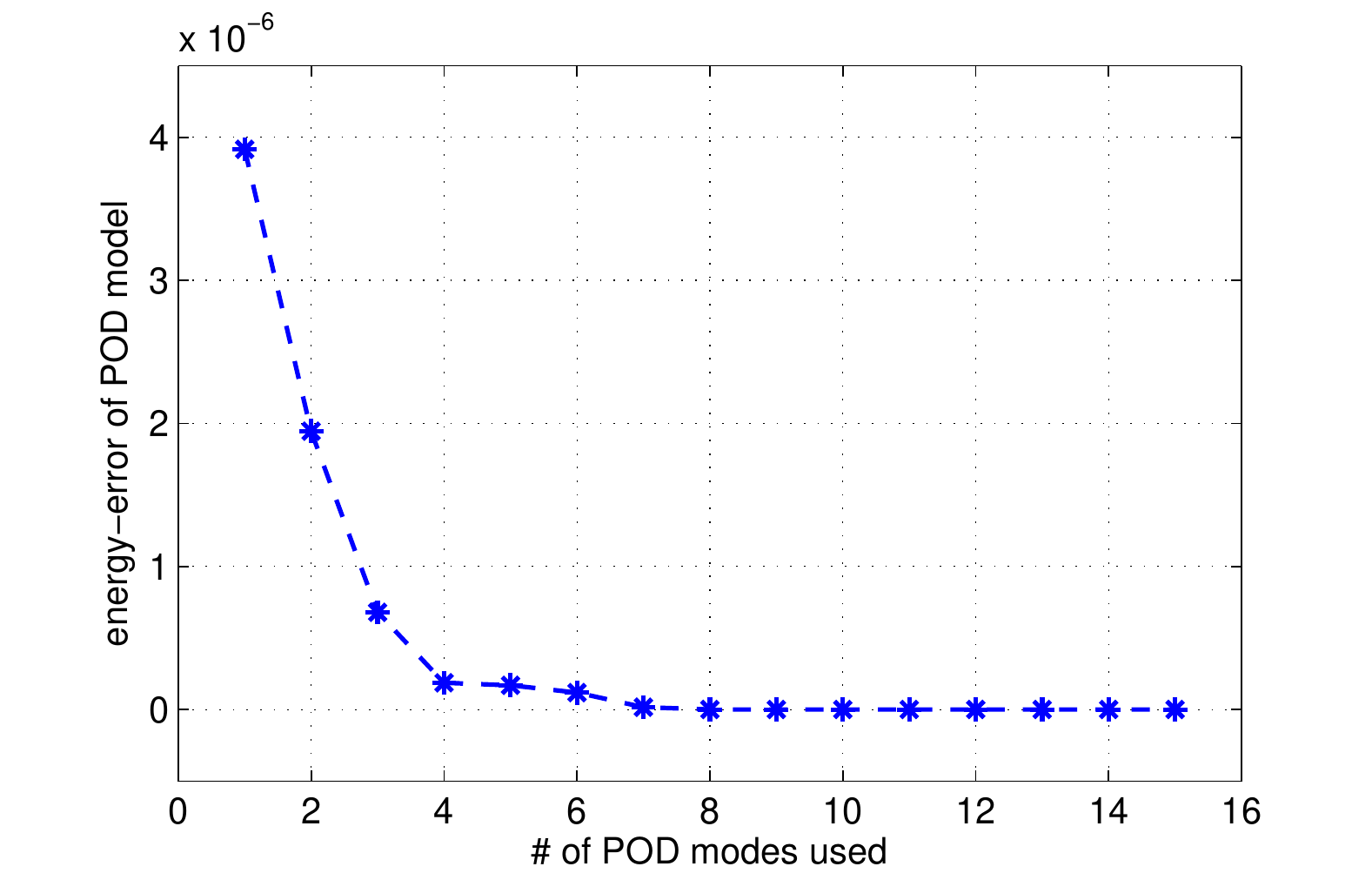}
\caption{1D NLS, Decay of the ROM errors :  solution (left),  Hamiltonian (right) }
\label{poderror}
\end{figure}

\begin{figure}[htb!]
\centering
\includegraphics[width=0.40\textwidth]{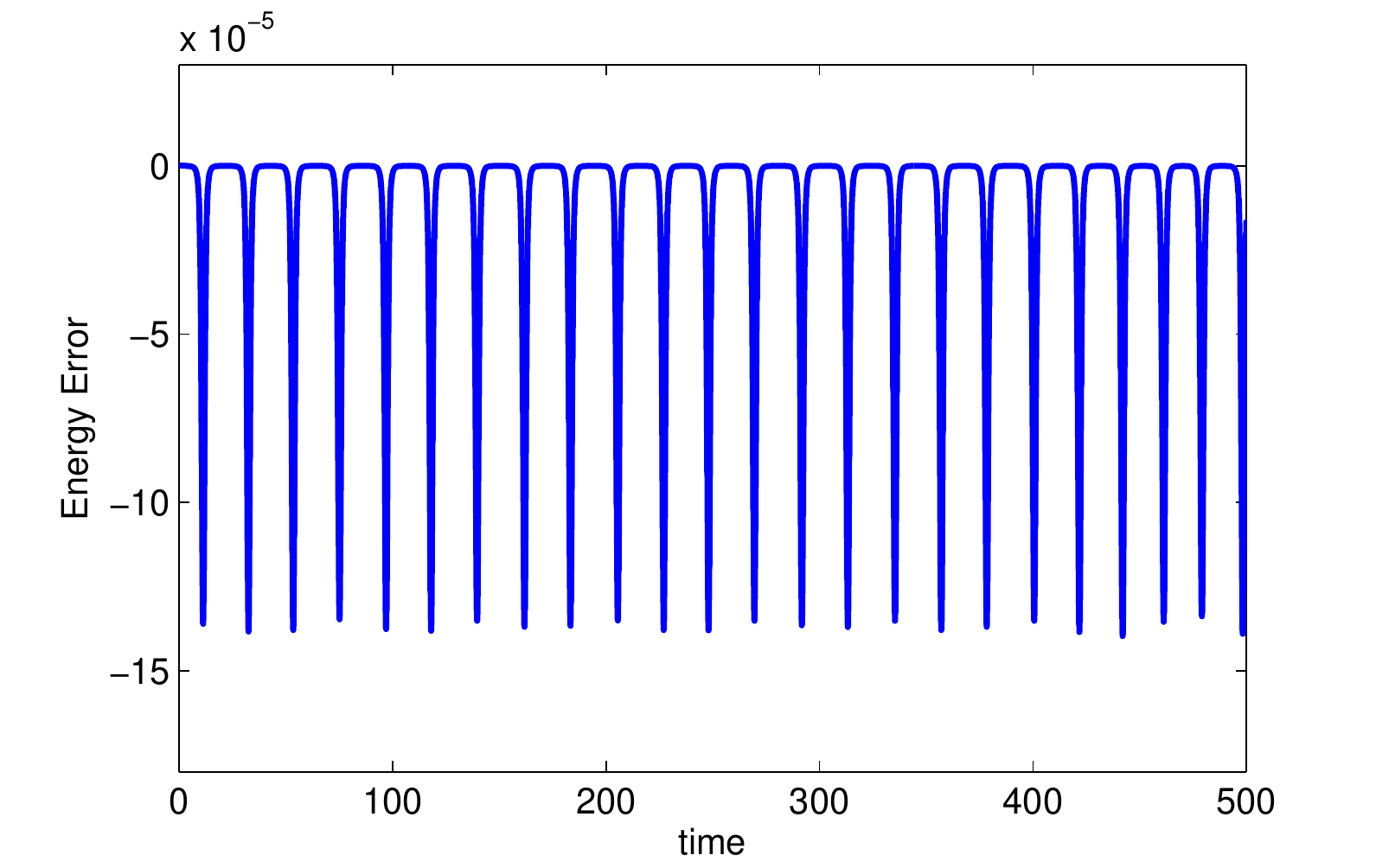}
\includegraphics[width=0.40\textwidth]{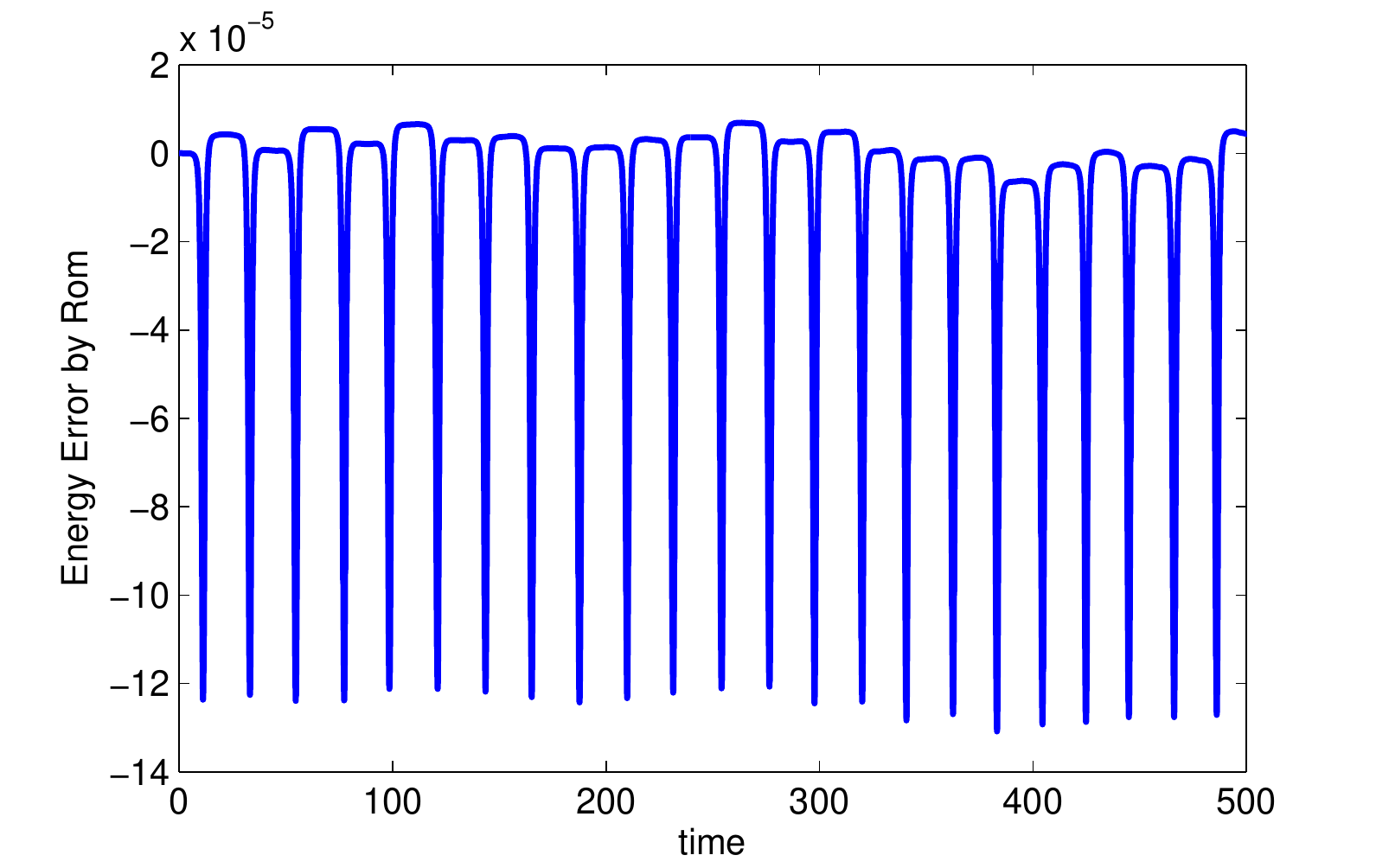}
\caption{1D NLS, Energy error: full-order model (left), ROM with 4 POD modes (right) }
\label{podeng}
\end{figure}

\begin{figure}[htb]
\centering
\includegraphics[width=0.40\textwidth]{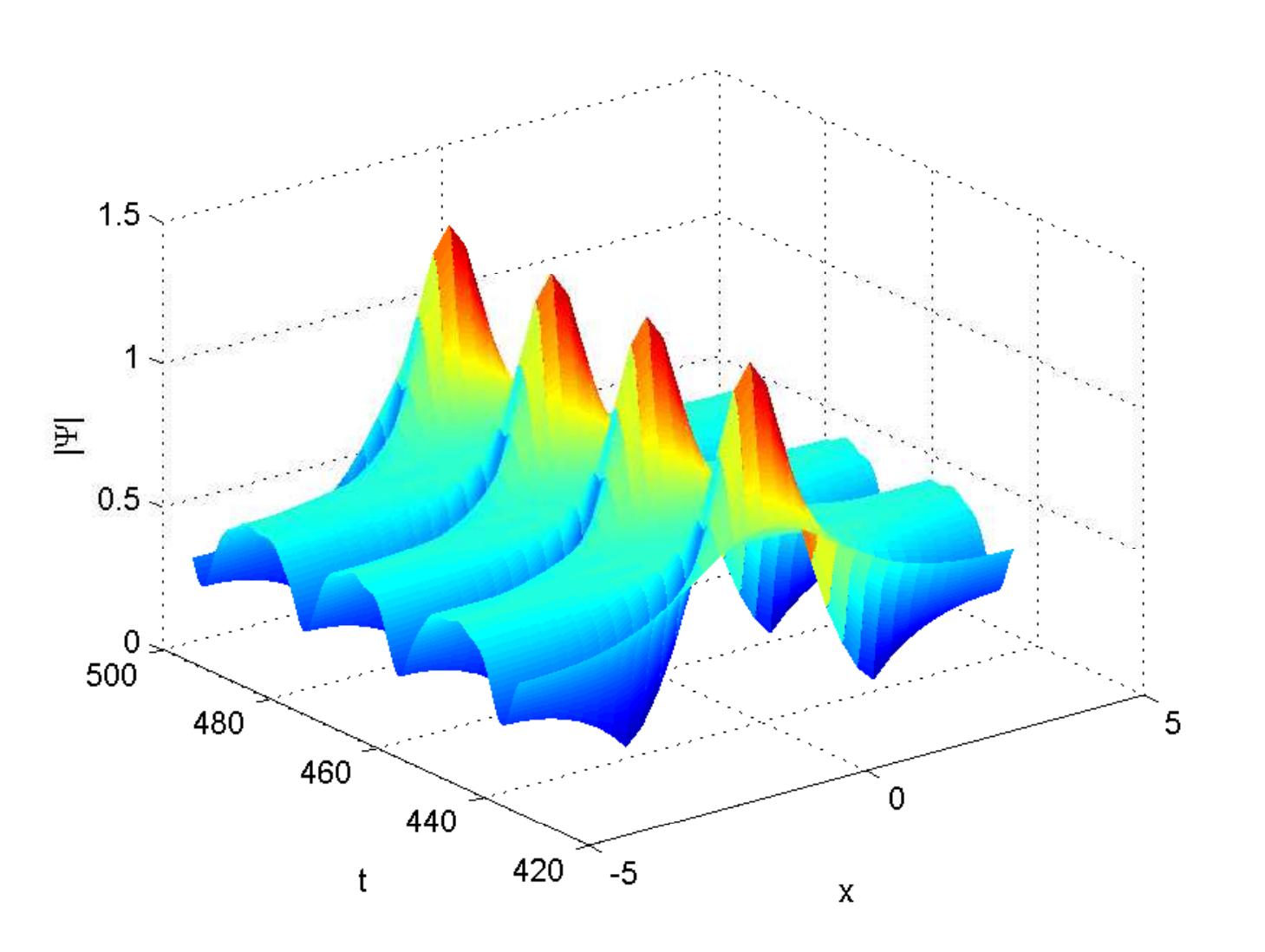}
\includegraphics[width=0.40\textwidth]{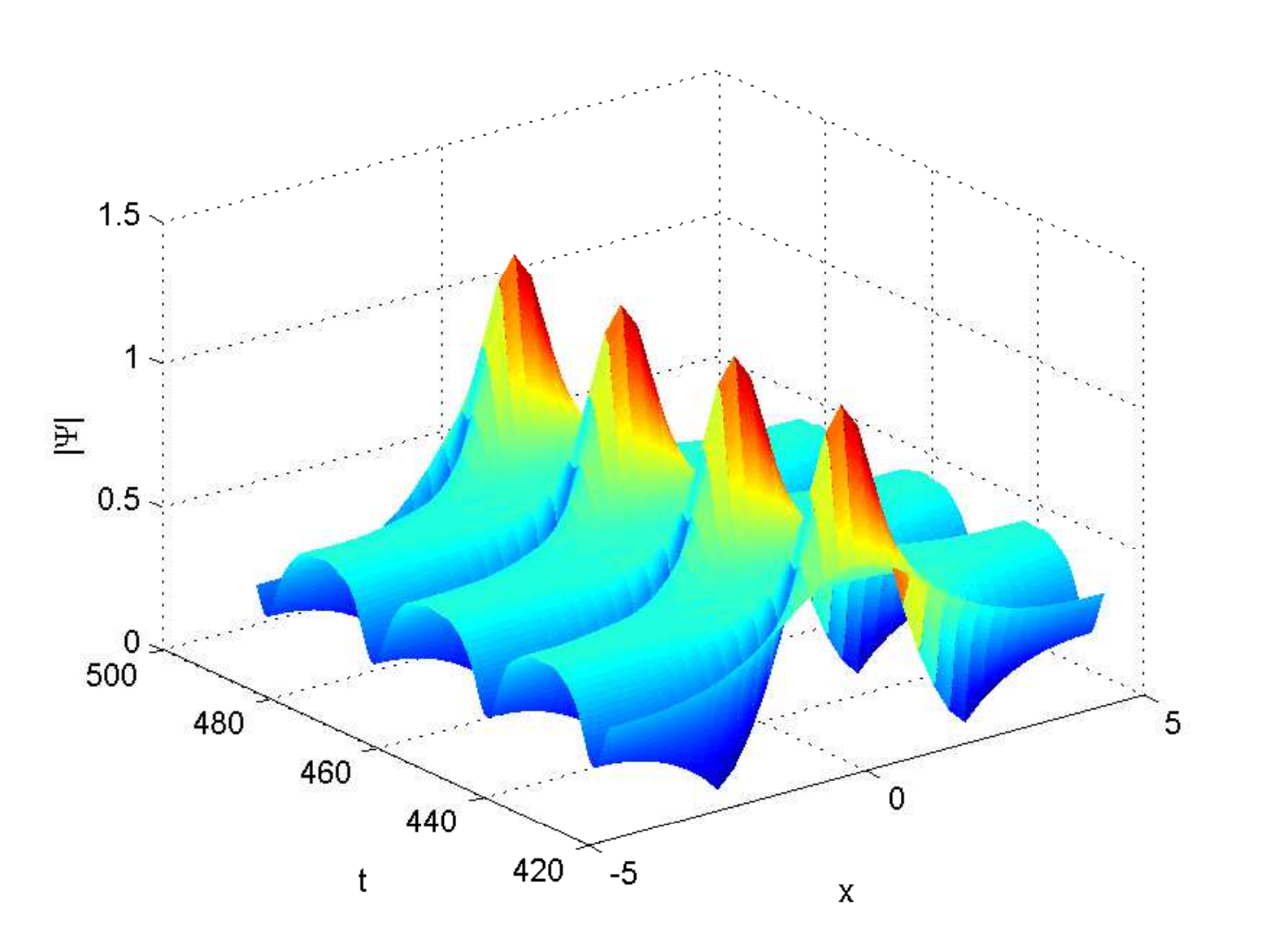}
\caption{1D NLS, Envelope of the approximate solution $\mid\psi\mid$: full-order model (left), ROM with 4 POD modes (right)}
\label{podsol}
\end{figure}

 For the one-dimensional NLS equation (\ref{nlsdenk}), we have taken the example in \cite{islas01gin} with $\gamma =2$
 and the periodic boundary conditions in the interval $[-L/2,L/2]$ with $L=2\sqrt{2}\pi$. The initial conditions are given
 as $p(x,0)=0.5(1+0.01\cos (2\pi x/L))$, $q(x,0)=0$. As mesh sizes in space and time,
 $\Delta x=L/32 $ and $\Delta t=0.01$ are used, respectively.
 Time steps are bounded by the stability condition  for the splitting method \cite{weidmann86sms};
  $\Delta t<\frac{2\Delta x^2}{L}$ where $L$ is the period of the problem. The discretized Hamiltonion is given by (\ref{hamilton}) with $\gamma =2$.\\
The singular values of the snapshot matrix are rapidly decaying (Figure \ref{svd}) so that  only few POD modes would  be sufficient to approximate the fully discretized NLS equation. In Figure \ref{poderror}, the  relative errors are plotted.
As expected with increasing number of POD basis functions $l$, the errors in the energy and
the errors between the discrete solutions of the fully discretized NLS equation and the
reduced order  model decreases which confirm the error analysis given in  Section 3.
In  Figure \ref{podeng} and \ref{podsol}, the evolution of the Hamiltonian error and the numerical solution at time $t=500$ are shown for the POD basis with $l=4$, where 99.99 \% of the energy of the system is well preserved. These figures confirm that the reduced model
 well preserves the Hamiltonian, and the numerical solution is close to the fully discrete solution.
\newpage
\subsection{Two-dimensional NLS equation}
\begin{figure}[htb!]
\centering
\includegraphics[width=0.4\textwidth]{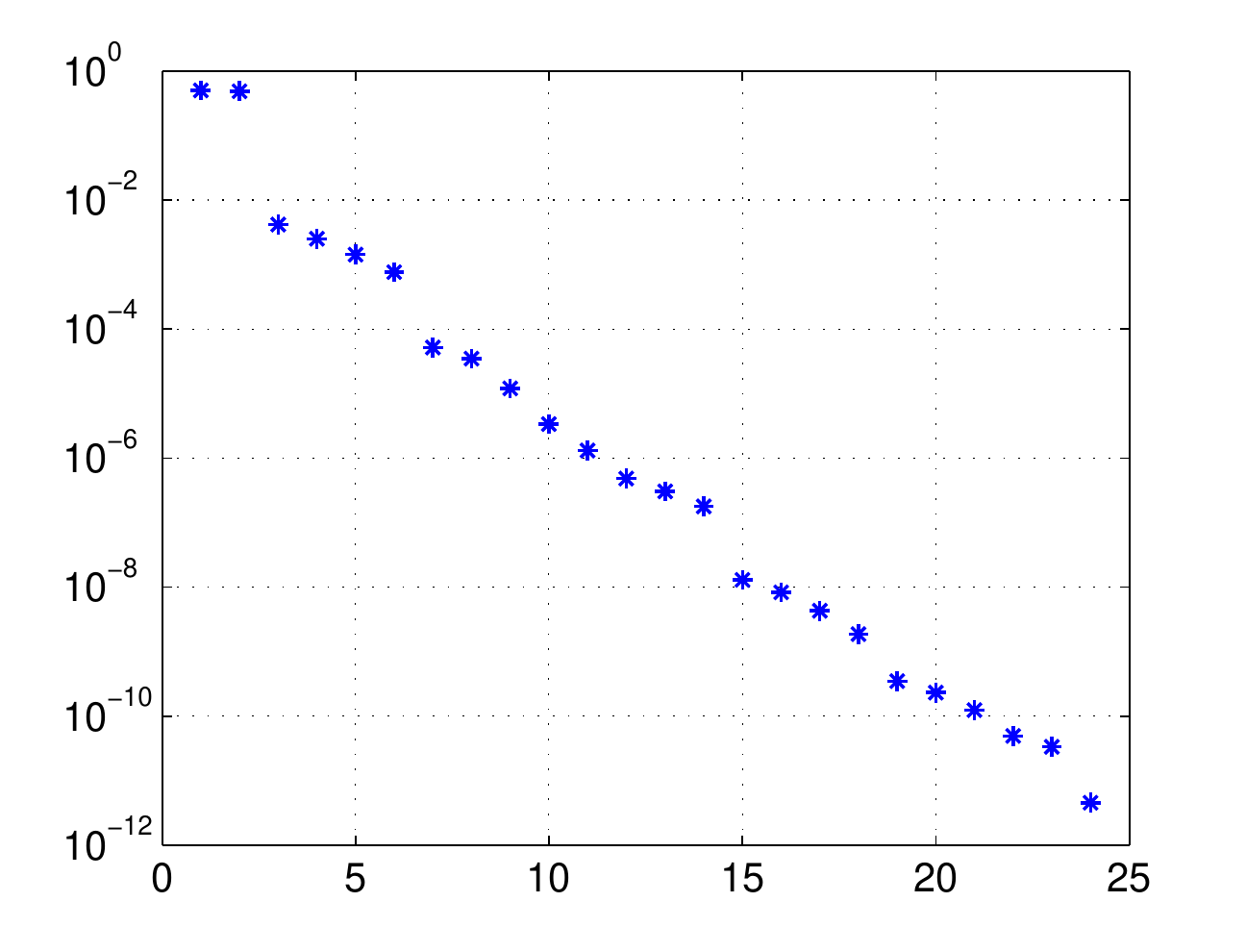}
\caption{2D NLS, Decay of the singular values}
\label{2dnlssv}
\end{figure}

\begin{figure}[htb!]
\centering
\includegraphics[width=0.40\textwidth]{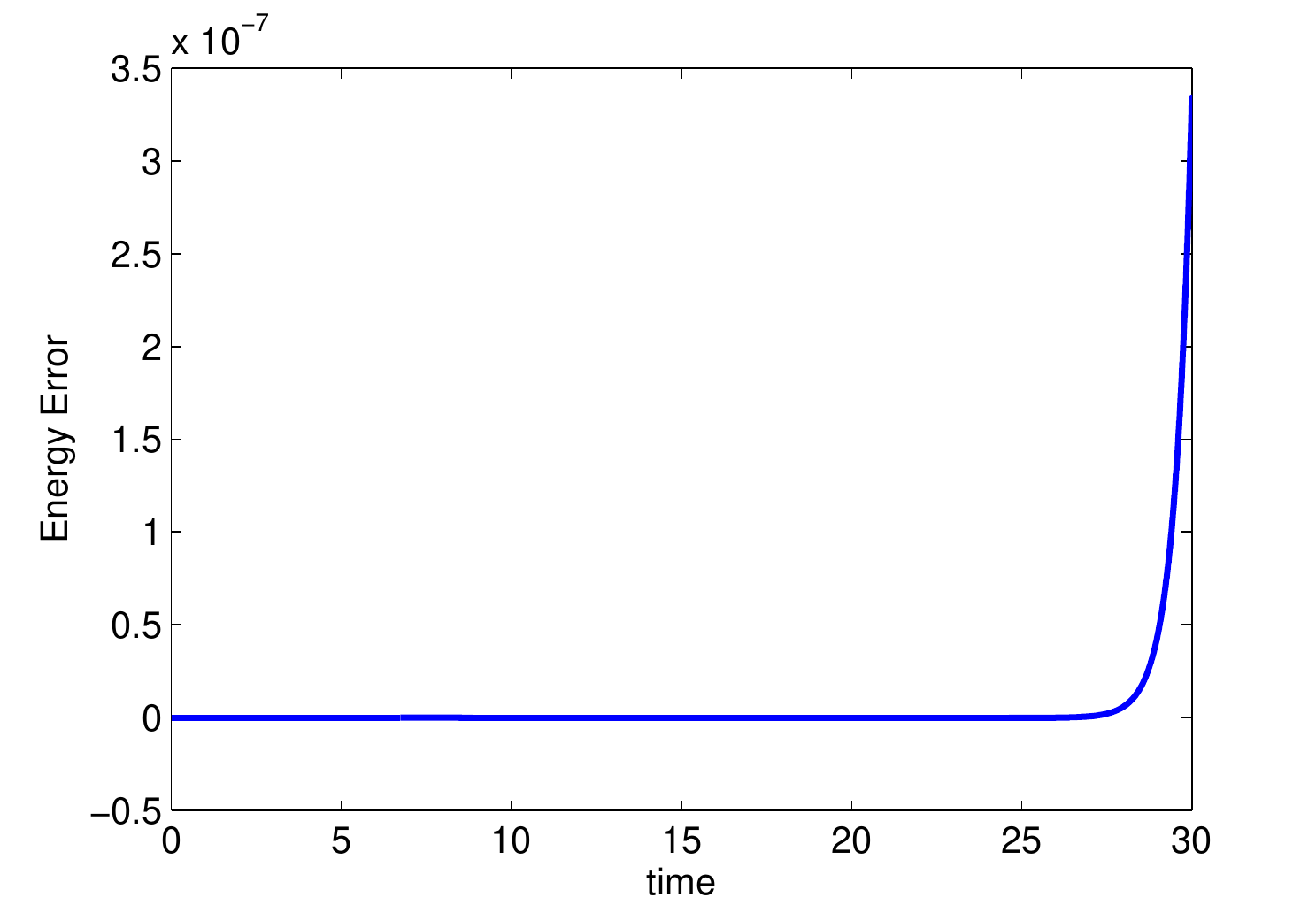}
\includegraphics[width=0.40\textwidth]{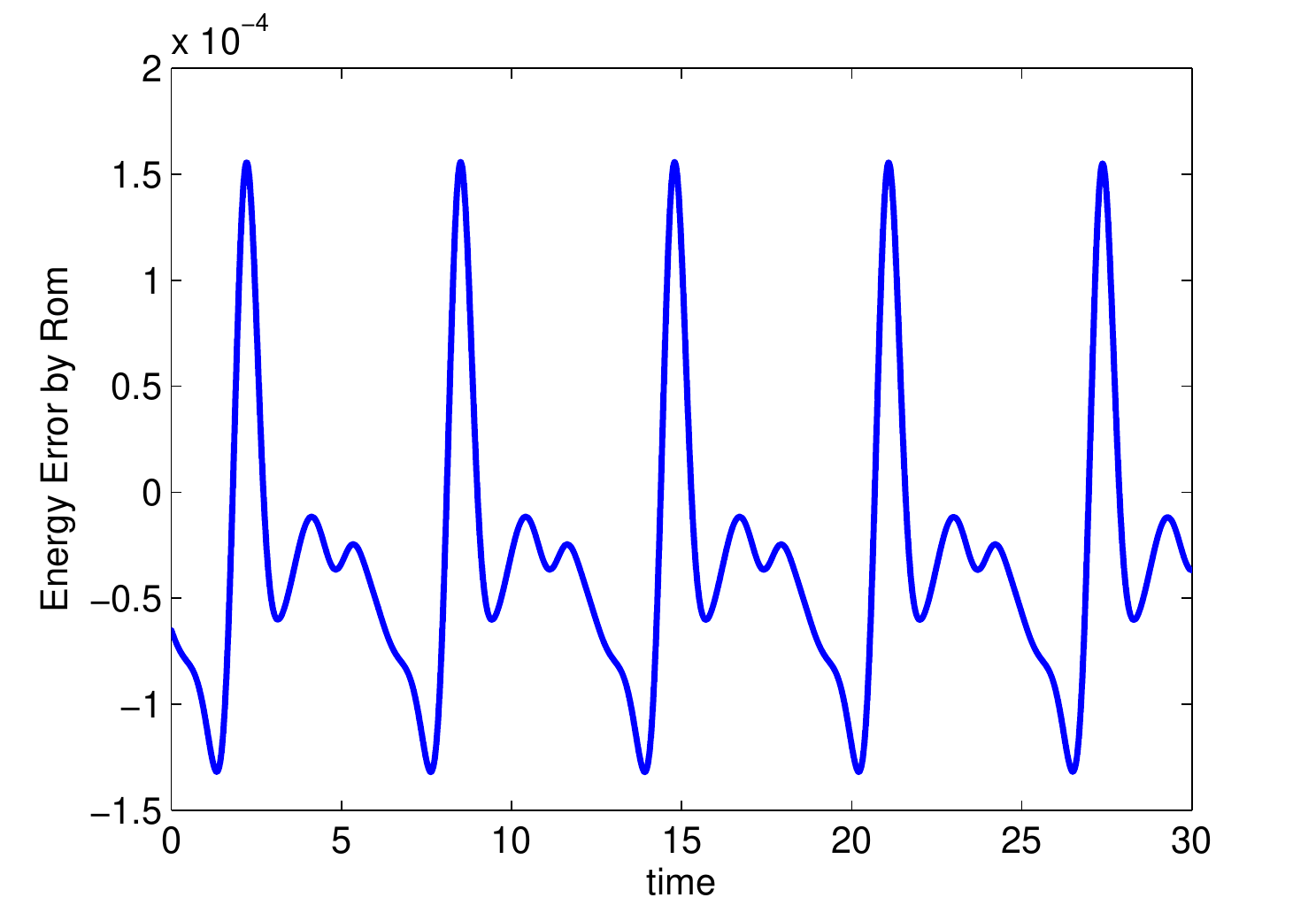}
\caption{2D NLS, Energy error: full-order model (left), ROM with 3 POD modes (right)}
\label{2d_error}
\end{figure}

We consider the following two-dimensional NLS equation \cite{chen11mss}
$$
\mathrm{i}\Psi_t+\Psi_{xx}+\Psi_{yy}+|\Psi|^2\Psi=0 \qquad \text{on } \;[0,2\pi]\times [0,2\pi]
$$
with the exact solution, $ \Psi(x,y,t)=\exp(\mathrm{i}(x+y-t))$.\\
The mesh size for  spatial discretization and time step size are taken as $\Delta x= \Delta y =2\pi /80 $ and $\Delta t=0.001$, respectively.
The discrete Hamiltonian is given by
\begin{eqnarray*}
H = &\Delta x\Delta y & \sum_{i,j=1}^m  \frac{1}{2} \left (\left ( \frac{p_{i+1,j} - p_{i,j}}{\Delta x} \right)^2 + \left ( \frac{q_{i+1,j} - q_{i,j}}{\Delta x} \right)^2  +\left ( \frac{p_{i,j+1} - p_{i,j}}{\Delta y} \right)^2 + \left ( \frac{q_{i,j+1} - q_{i,j}}{\Delta y}\right)^2 \right)\\
& & - \frac{1}{4} \left( p_{i,j}^2 + q_{i,j}^2 \right)^2
\end{eqnarray*}
Only 3 POD modes were sufficient to capture almost all of the energy of the system (Table \ref{2d_table}). A comparison  of the Hamiltonian errors in long term computation shows that the reduced order model with a few POD modes preserve the energy of the system very well (Figure \ref{2d_error}). The singular values of 2D NLS are decreasing not continuously as for 1D NLS equation (Figure \ref{2dnlssv}).
\begin{table}[htb!]
\centering
\begin{tabular}{|c|c|c|c|}
\hline
\# POD & Info (\% ) &  ROM Hamiltonian error & ROM error \\
\hline \hline
1 & 51.65    &    8.181e-002  &  2.770e+001  \\
2 & 99.995   &    6.116e-007  &  1.040e-003  \\
3 & 99.998   &    4.164e-007  &  1.134e-003 \\
\hline
\end{tabular}
\caption{2D NLS, RIC and errors for the real part of the solution}
\label{2d_table}
\end{table}
\newpage
\subsection{Coupled NLS equation }
\begin{figure}[htb!]
\centering
\includegraphics[width=0.40\textwidth]{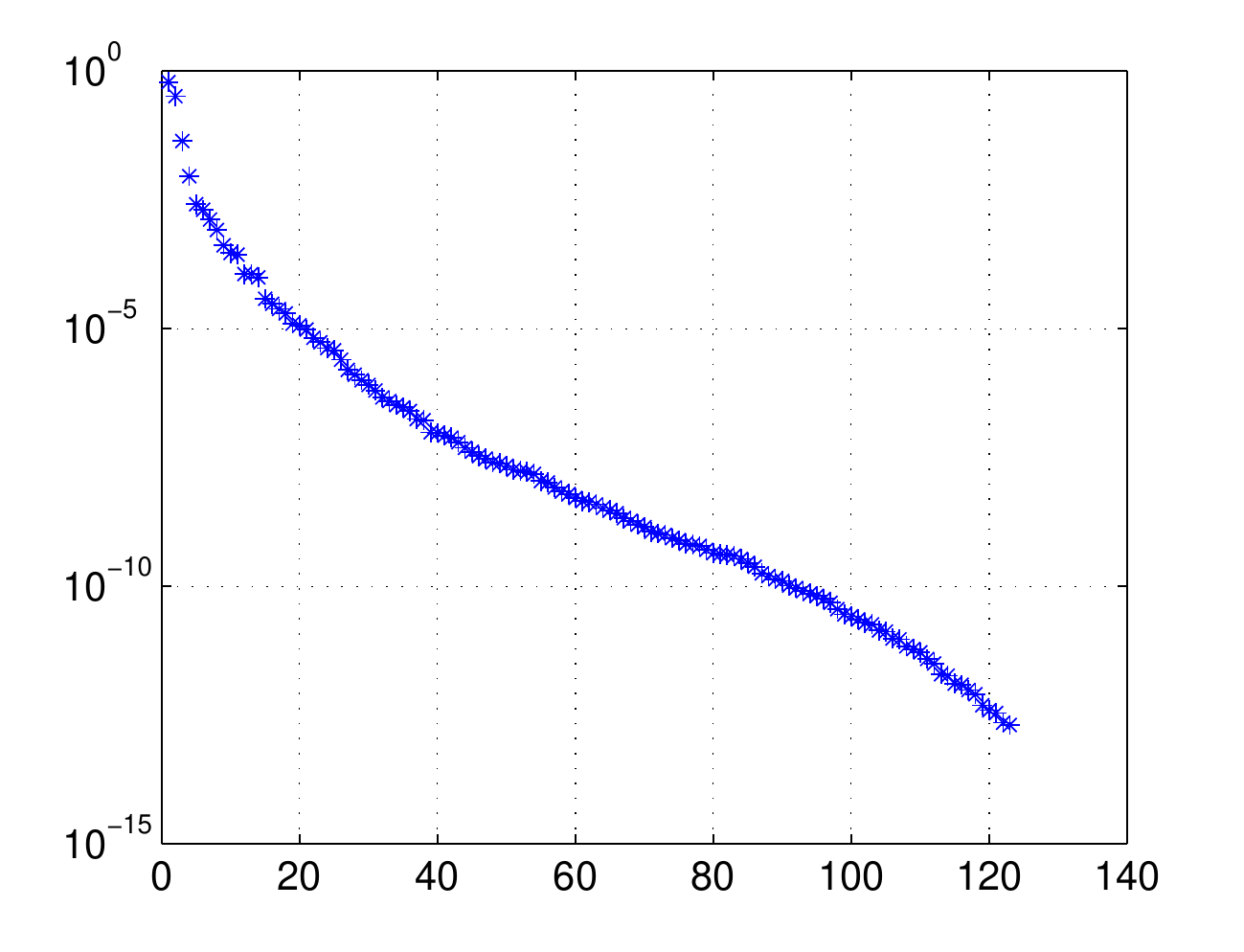}
\caption{Coupled NLS, Decay of the singular values for the real part of $\Psi_1$}
\label{cnlssv}
\end{figure}

\begin{figure}[htb!]
\centering
\includegraphics[width=0.40\textwidth]{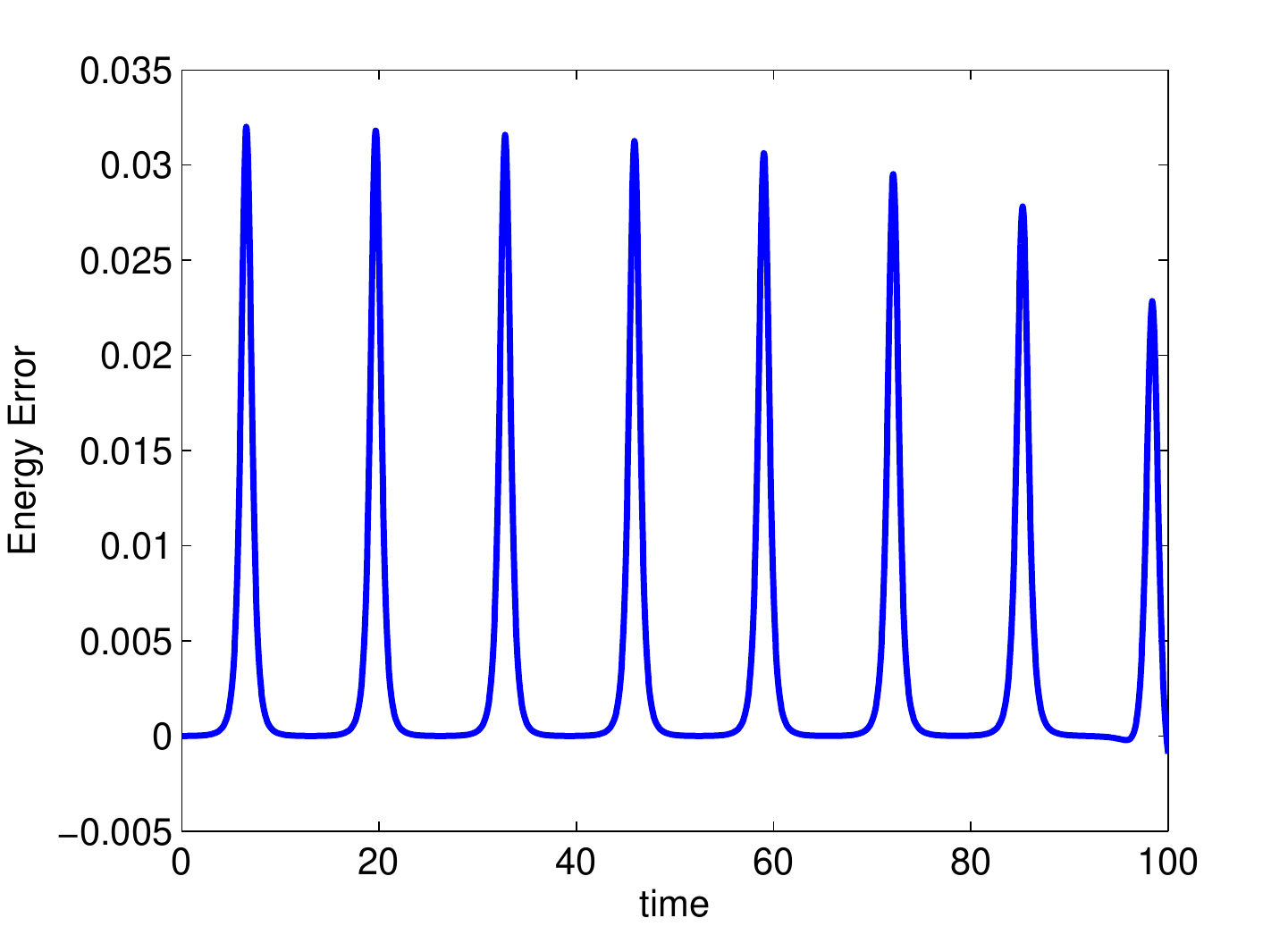}
\includegraphics[width=0.40\textwidth]{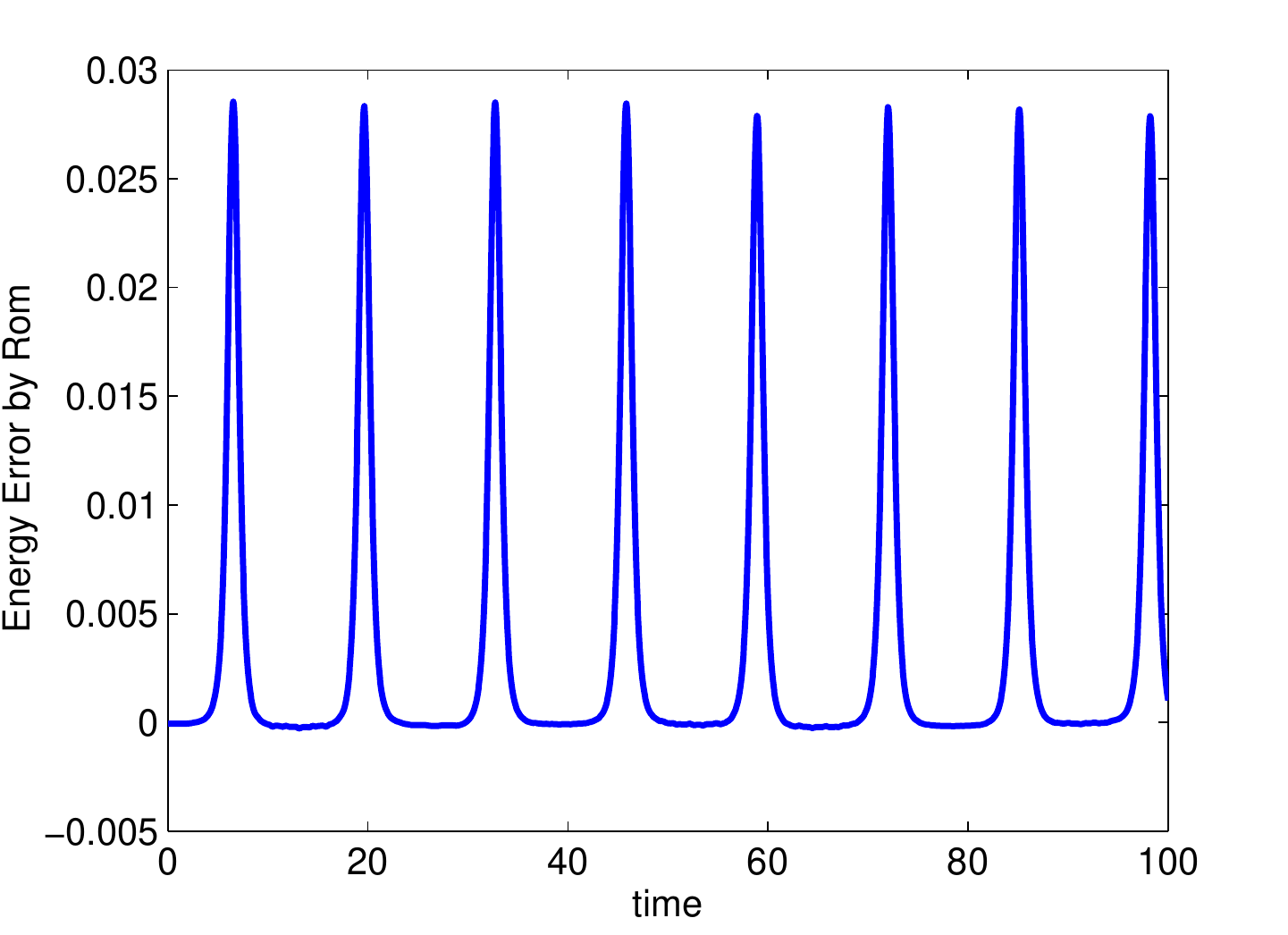}
\caption{Coupled NLS, Hamiltonian error: full-order model (left) and ROM with 5 POD modes (right)}
\label{cnls_ham}
\end{figure}

\begin{figure}[htb!]
\centering
\includegraphics[width=0.40\textwidth]{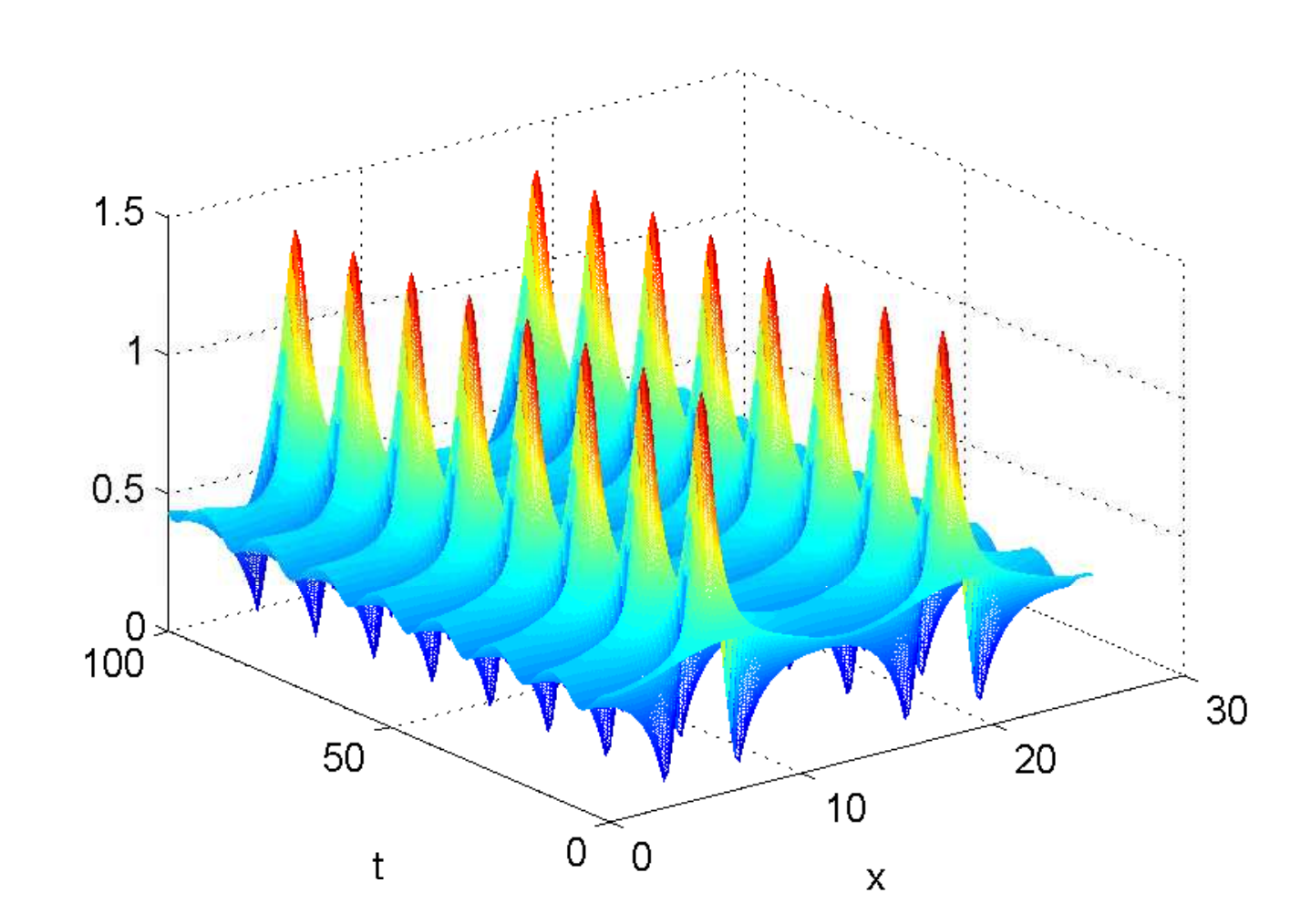}
\includegraphics[width=0.40\textwidth]{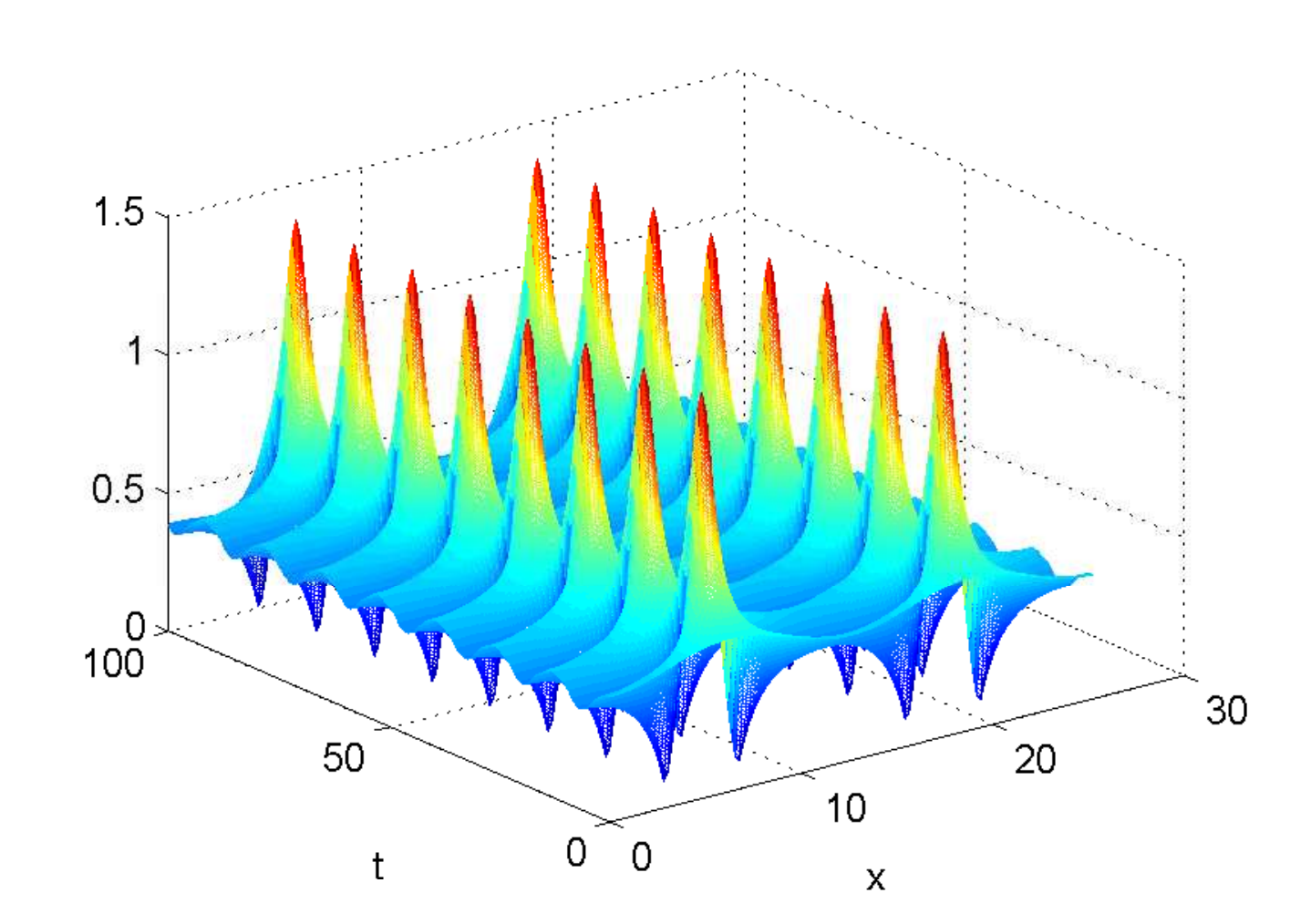}
\caption{Coupled NLS, Interaction of solitons $|\Psi_1|$ and $|\Psi_2|$: full-order model (left) and ROM with 5 POD modes (right)}
\label{cnls_sol}
\end{figure}

We consider two coupled NLS equations (CNLS) with elliptic polarization with plane wave solutions \cite{ayhan07smm}
\begin{equation}
\mathrm{i}\frac{\partial\Psi_1}{\partial t}+\frac{\partial^2\Psi_1}{\partial x^2}+(|\Psi_1|^2+|\Psi_2|^2)\Psi_1=0,\quad
\mathrm{i}\frac{\partial\Psi_2}{\partial t}+\frac{\partial^2\Psi_2}{\partial x^2}+(|\Psi_2|^2+|\Psi_1|^2)\Psi_2=0
\end{equation}
using the initial conditions
$$
\Psi_1(x,0)=(0.5)(1-0.1\cos (0.5x)), \qquad \Psi_2(x,0)=(0.5)(1-0.1\cos (0.5x)).
$$
The equations  are solved over the space $[0,8\pi]$ and time interval $[0,100]$, respectively, with the mesh size and time steps
$
dx=8\pi /128, \; \Delta t=0.05.
$
The discrete Hamiltonian is given as \cite{ayhan07smm}
\begin{eqnarray*}
H = &\Delta x & \sum_{j=1}^m  -\frac{1}{2} \left (\left ( \frac{p^1_{j+1} - p^1_j}{\Delta x} \right)^2 +  \left ( \frac{q^1_{j+1} - q^1_j}{\Delta x} \right)^2
+ \left ( \frac{p^2_{j+1} - p^2_j}{\Delta x} \right)^2 + \left ( \frac{q^2_{j+1} - q^2_j}{\Delta x} \right)^2 \right) \\
& & + \frac{1}{4} \left( ( (p^1_j)^2 + (p^2_j)^2)^2 + ( (q^1_j)^2 + (q^2_j)^2)^2 \right) + \frac{1}{2}\left(  ((p^1_j)^2 + (p^2_j)^2)((q^1_j)^2 + (q^2_j)^2) \right)
\end{eqnarray*}
where $p^1, q^1$ and $p^2, q^2$ denote the real and imaginary parts of $\psi_1$ and $\psi_2$, respectively.\\
Figure \ref{cnls_ham} \& \ref{cnls_sol} and Table  \ref{cnlstab}  show that only few POD modes are necessary to capture the dynamics of the CNLS equation.
The singular values are  decreasing not so rapidly (Figure \ref{cnlssv}) as in case of single 1D and 2D  NLS equations.
\begin{table}[htb!]
\centering
\begin{tabular}{|c|c|c|c|}
\hline
\#POD & RIC(\%) & ROM Hamilton error & ROM error\\
\hline \hline
 2  & 99.58 &  1.879e-004  &  5.060e-001  \\
 3  & 99.98 &  1.865e-004 &  3.761e-001\\
 4  & 99.99 &  1.213e-004  &  6.491e-002  \\
 5  & 99.99 &  2.825e-005  &  3.919e-003\\
 \hline
\end{tabular}
\caption{Coupled NLS, RIC and errors for the real part of $\Psi_1$}
\label{cnlstab}
\end{table}
\section{Conclusions}
A reduced model is derived for the NLS equation by preserving the Hamiltonian structure. A priori error estimates are obtained for the mid-point rule as time integrator for the reduced dynamical system. Numerical results show that the energy and the phase space structure of the three different NLS equations are well preserved by using few POD modes. The number of the POD modes containing most of the energy depends on the decay of the singular values of the snapshot matrix, reflecting the dynamics of the underlying systems. In a future work, we will investigate the dependence of the ROM solutions on parameters for the CNLS equation by performing a sensitivity analysis.

\end{document}